\documentclass[12pt]{amsart}
\usepackage{amsmath,amsfonts,amssymb,amsthm,amstext}
\usepackage{geometry}
\usepackage{hyperref}
 \usepackage{amsfonts}
\usepackage{amsmath}
\usepackage{amssymb}
\usepackage{amstext}
\usepackage{amsthm}  
\usepackage{epigraph}    
\newskip\stdskip                      
\theoremstyle{plain}
\newtheorem{Theorem}{Theorem}[section]

\newtheorem{Corollary}[Theorem]{Corollary}

\newtheorem{Question}{Question}
\theoremstyle{definition}
\newtheorem{Definition}[Theorem]{Definition}
\newtheorem{Remark}[Theorem]{Remark}
\newtheorem{Example}[Theorem]{Example}

\newcommand{\C}{\mathbb C}
\newcommand{\T}{\mathbb T}

\newcommand{\bS}{\mathbb S}

\begin{document}
%--------------------------------------------------------------------%
% Article 1_park
%--------------------------------------------------------------------%
 \title{Some problems on intrinsically harmonic forms}
\author{Gianluca Bande}
\address{Dipartimento di Matematica e Informatica, Universit\`a degli Studi di Cagliari, Cagliari, Italy}
\email{gbande@unica.it}

\thanks{The author is member of the Italian National Group G.N.S.A.G.A. of INdAM.The  author is partially supported by the project DAMPAI - Fondo Crescita Sostenibile - B29J24001610005 - F/350361/01-02/X60}

\date{\today}

\maketitle

   \begin{center}
 {\small  To the memory of Francesco Mercuri}
  
    \end{center}

\begin{abstract}
  In this short note we recall the definition of intrinsically harmonic forms, some known results and some open problems.
\end{abstract}

%  \classification{primary 57R22; secondary 53C65, 53C05, 53C15}

\section{Introduction}

A differential form $\omega$ on an $n$-dimensional manifold $M$ is called {\it intrinsically harmonic} \cite{C} if it is harmonic with respect to some Riemannian metric $g$. An intrinsically harmonic form is a fortiori closed and then the main problem is to give necessary and sufficient conditions under which a closed form is intrinsically harmonic.

In fact only the forms of degrees $1$ and $n-1$ have been quite well understood. A  classical theorem of Calabi \cite{C} answers the question for $1$-forms with non-degenerate zeros and Honda \cite{H} proved the dual case of $(n-1)$-forms.

In 2007 Volkov \cite{V} gave a complete characterization for closed $1$-forms, but in general the forms of degrees strictly between $1$ and $n-1$ present additional problems. 

In this short survey I recall some known results and a few problems concerning in particular the case of $2$-forms. In what follows we will assume that the manifolds we are woking with are closed, connected and $n$-dimensional.

%%%%%%%%%%%%%%%%%%%%%

\section{Intrinsically harmonic forms}\label{sec: preliminari}

The problem of giving an intrinsic characterization of harmonic forms has been introduced by Calabi in \cite{C}. The case of $0$-forms is quite trivial because we are working with closed manifold, so only constant function are intrinsically harmonic. On the other side, every volume form is intrinsically harmonic because it is closed and can be realized as a Riemannian volume form. The first interesting case arises for $1$-forms.

Calabi himself proved the following characterization of {\it generic} $1$-forms (a $k$-form is said to be generic when, as section of the vector bundle of $k$-forms, is transverse to the zero section):
\begin{Theorem}\label{th: calabi}
A generic closed $1$-form $\alpha$ is intrinsically harmonic if and only if the following conditions are satisfied:
\begin{enumerate}
\item $\alpha$ is locally intrinsically harmonic (i.e. intrinsically harmonic on a neighborhood of the set of its singularities),
\item  for every non-singular point $p$ there exists an embedded circle passing through $p$ and transverse to the hyperplane determined by $\alpha$ ($\alpha$ is said transitive).
\end{enumerate}
\end{Theorem}
Calabi also showed that a nowhere vanishing closed $1$-form is necessarily transitive and then intrinsically harmonic.

In a similar fashion, Honda \cite{H} proved the dual case of $(n-1)$-forms, where transitivity here is understood as the existence, on non-singular points, of a codimension $1$ closed submanifold transverse to the line field determined by the $(n-1)$-form. 

Mercuri and Fran\c ca \cite{EM} proved that transitive nowhere vanishing $(n-1)$-forms are intrinsically harmonic.

\begin{Remark}\label{rem: n-1 forms}
In contrast with the case of closed $1$-forms, for never-vanishing closed $(n-1)$-forms it is not true that they are automatically transitive and then intrinsically harmonic. A simple example is given by the pullback in $\bS^3$ of the volume form of $\bS^2$ via the projection of the Hopf fibration (see \cite{EF}, for example).
\end{Remark}

In 2007, Volkov \cite{V} was able to drop the genericity condition used by Calabi and gave a complete characterization for closed $1$-forms. 

Volkov also pointed out that, in general, the forms of degrees strictly between $1$ and $n-1$ present additional problems. One of this difficulties is illustrated by an interesting example which we recall in Section \ref{sec: IH2forms}.

\section{Intrinsically harmonic two-forms and so on}\label{sec: IH2forms}

For forms with degrees strictly between $1$ and $n-1$, non-degeneracy means that the forms have no zeros at all (just for dimensions reasons and transversality).

Simple examples of generic $2$-forms are given by any symplectic form which is harmonic with respect to any compatible metric. In fact, its Hodge-dual is, up to a constant, a power of the symplectic form itself. 

Another simple example is given by any of the two forms composing a {\it symplectic pair} defined as follows:

\begin{Definition}[\cite{BGK, BK}]
  Let $M$ be a $2m$-dimensional manifold. A pair of {\em closed} $2$-forms 
$(\omega,  \eta)$ is called a {\it symplectic pair} of type $(k, m-k)$ (for $0<k<m$) if they have constant 
ranks $2k$ and $2(m-k)$ respectively, and moreover $\omega^{2k} \wedge \eta^{2(m-k)}$ 
is a volume form.
\end{Definition}
%A symplectic pair gives rise to two symplectic forms
%$$\Omega_+=\omega + \eta, \quad \Omega_-= \omega - \eta$$ 
%on $M$ and on $(-1)^{n-p}M$ respectively, where $-M$ denotes the oriented manifold 
%obtained by reversing the orientation of $M$.  
%To make the definition interesting, we will assume that $k>0$ and $n>k$. Then, when $M$ 
%has dimension four  --- the case of interest in the present article --- a symplectic pair on $M$ 
%can only be of type $(1,1)$ and, in particular, $M$ is symplectic for both orientations. 

%In dimension four, a symplectic pair $(\omega, \eta)$ can be equivalently 
%defined by a pair of symplectic forms $(\Omega_+,  \Omega_-)$ 
%satisfying 
%$$
%\Omega_+^2=-\Omega_-^2 \quad , \quad \Omega_+\wedge \Omega_-=0 \; .
%$$
%The symplectic pair is then given by $$\Big(\frac{\Omega_+ +\Omega_-}{2} ,
%\frac{\Omega_+ -\Omega_-}{2} \Big) \, ,$$
%and we will say that $(\Omega_+,  \Omega_-)$ {\it arises from a symplectic pair}.
%\begin{rem}\label{nondefinite}
%If $M$ is a closed four-manifold which is symplectic for both orientations, then 
%$b_\pm(M)>0$ by Equation~\eqref{eq: two forms}. In particular, $\C \PP^2$ does 
%not admit a symplectic pair. 
%\end{rem} 

The kernels of $\omega$ and $\eta$ are integrable complementary distributions and therefore
integrate to a pair of transverse foliations ${\mathcal F}_\omega$ and ${\mathcal F}_\eta$, called {\it characteristic foliations}.
%such that
%$$T{\mathcal F}_\omega=\ker \omega \quad \text{and} \quad T{\mathcal
%  F}_\eta =\ker \eta.$$ 
%Each form is symplectic on the leaves of the foliation induced 
%by the other form and moreover ${\mathcal F}_\omega$ and ${\mathcal F}_\eta$ are 
%symplectically orthogonal with respect to both the symplectic forms $\Omega_{+}$ and $\Omega_{-}$.
%
%By Rummler and Sullivan's criterion (see \cite{G}), the characteristic foliation of a closed $2$-form is taut, which means that there exists a Riemanniann metric for which the leaves are minimal. 

For a symplectic pair it is possible to construct so-called {\it compatible} Riemannian metrics. With respect to these metrics, the foliations are orthogonal, and both have minimal leaves. Moreover, as for a symplectic form, each form of the pair is harmonic with respect to a compatible metric.

%%%%%%%%%%%
%%%%%%%%%%%

Concerning the problem of characterizing generic $2$-forms, the first interesting case comes up in dimension four, and at a first glance one could consider only forms of constant rank.

On a $4$-dimensional manifold, a non-vanishing $2$-form of constant rank has either rank $2$ or $4$. Zero forms and symplectic forms (of rank $4$) are intrinsically harmonic, so let us consider $2$-forms of rank $2$. 

Concerning closed forms of rank $2$, Volkov gave the following example (see  \cite[Example 1]{V}):
\begin{Example}\label{ex:volkov}
Let $\C P^1\tilde\times \C P^1$ be the non-trivial $\C P^1$-bundle over $\C P^1$. On the base space of $\C P^1\tilde\times \C P^1$, consider a volume form  and pull it back on the total space. Let's call $\Omega$ the pull back. Of course $\Omega$ is closed, has constant rank $2$ and its kernel foliation is given by the fibers of the bundle. Nevertheless $\Omega$ cannot be intrinsically harmonic because otherwise its Hodge dual would also have constant rank two, and would determine a foliation transverse to the fibers and therefore a foliated bundle. Foliated bundles are classified (up to conjugation) by a representation of the fundamental group of the base space into the diffeomorphism group of the fiber, but the simply connectedness of the base forces the bundle to be trivial.
\end{Example}

The discussion after Example 1 in \cite{V}, suggested in some sense the following result:
\begin{Theorem}[\cite {EM}]\label{th Elizeu} Let M be a manifold with a closed $k$-form $\omega$ of rank $k$.
Then $\omega$ is intrinsically harmonic if and only if there exists a closed form $\eta$ of rank
$n-k$ such that $\ker \omega \cap \ker \eta= \{0\}$.
\end{Theorem}

If we restrict to $2$-forms on $4$-manifolds then we have:
\begin{Theorem}[\cite{B}]\label{th pr 2}
A closed $2$-form of constant rank $2$ on an orientable closed four dimensional manifold $M$ is intrinsically harmonic if and only if it is part of a symplectic pair.
\end{Theorem}

Still working with bundles, Fran\c ca and Mercuri were able to characterize flat principal bundles in term of intrinsic harmonicity:
\begin{Theorem}[\cite {EM}]\label{th EM -constant rank} Let $B= \{B,p,M\}$ be a diﬀerentiable principal circle bundle with M closed and orientable. Then $B$ admits a flat connection if and only if the pullback of a volume form of $B$ is intrinsically harmonic.
\end{Theorem}
%%%%%%%%%
%%%%%%%%%
As an application of the theory of intrinsically harmonic forms, we close this section with the following result which characterizes the $n$-torus among closed manifolds supporting $(n-1)$ linearly independent closed $1$-forms:
\begin{Theorem}\cite{EF}
Let $M$ be a closed $n$-dimensional manifold supporting $(n-1)$ linearly independent closed $1$-forms $\lambda_1,\cdots,\lambda_{n-1}$. If $\omega:= \lambda_1 \wedge \cdots \wedge \lambda_{n-1}$ determines a non-zero cohomological class, then it is intrinsically harmonic. In particular, there exists a closed $1$-form $\eta$ such that the set
$\{ \lambda_1,\cdots,\lambda_{n-1}, \eta\}$ is linearly independent, and therefore M is diﬀeomorphic to $\T^n$.
\end{Theorem}
%%%%%%%%
Concerning the link with vector fields, I recall the following result of Simi\'c (a different proof is also given in \cite{F}), related to the existence of global cross-sections for non-singular flows:
\begin{Theorem}\cite[Main Theorem]{SI}
Let $\Phi$ be a non-singular smooth flow on a smooth, compact, connected manifold $M$. Denote the infinitesimal generator of $\Phi$ by $X$ and assume that $\Phi$ preserves a smooth volume form $\Omega$. Then $\Phi$ admits a smooth global cross-section if and only $i_X\Omega$ is intrinsically harmonic.
\end {Theorem}
As a corollary one has the following:
\begin{Corollary}\cite[Corollary 1.5]{SI}\label{cor-simic}
Let $\Phi$ be a non-singular smooth flow on a smooth, compact, connected manifold $M$, with infinitesimal generator $X$. Assume that $\Phi$ preserves a smooth volume form $\Omega$. If $i_X \Omega$ is intrinsically harmonic, then $[i_X \Omega] \neq0 \in H^{n-1}_{deRham} (M)$.
\end{Corollary}

%%%%%%%%
%%%%%%%%
\section{Open problems}

The interest of Theorem \ref{th pr 2} is due to the fact that there are topological obstructions to the existence of a symplectic pair on closed manifolds. For example, the existence of a symplectic pair on a manifold implies that its second Betti number $b_2$ satisfies $b_2\geq 2$. Therefore $\C P^2$, which is symplectic, admits intrinsically harmonic $2$-forms of constant rank $4$ but no intrinsically harmonic $2$-forms of constant rank $2$.

%But on $\C P^2$ no closed $2$-form of constant rank $2$ exists, because there is no splitting of the tangent bundle into two rank $2$ subbundle.

A more interesting example is given by $\C P^2 \# \overline{\C P^2}$, the non-trivial $\C P^1$-bundle over $\C P^1$ of Example \ref{ex:volkov} (see \cite{MS} for an explanation). $\C P^2 \# \overline{\C P^2}$ fulfills all the basic topological obstructions to the existence of a symplectic pair but, by \cite{BG}, it admits no symplectic pair and thus we have:

\begin{Corollary}
$\C P^2 \# \overline{\C P^2}$ admits no intrinsically harmonic $2$-form of constant rank $2$. 
\end{Corollary}

Since $\C P^2 \# \overline{\C P^2}$ is symplectic, the only intrinsically harmonic $2$-forms of constant rank, can have rank $4$ or $0$. 

%An interesting problem is then to find an example of intrinsically harmonic $2$-form of non-constant rank $k$ with  $k\leq 2$ or $k\leq 4$ everywhere.  

The following natural question is answered by the positive in \cite{B}:
\begin{Question} Does $\C P^2 \# \overline{\C P^2}$ admit an intrinsically harmonic $2$-form of non-constant rank which is not symplectic and has at least rank $2$ in a point? 
\end{Question}
%The answer is yes and an example can be constructed as follows.
%
%\begin{Example}
%Let $\omega$ be  the pullback to $\C P^2 \# \overline{\C P^2}$ of the Fubini-Study volume form on $\C P^1$, and fix any Riemannian metric $g$ on $\C P^2 \# \overline{\C P^2}$. By the Hodge Theorem, the cohomology class of $\omega$ has a unique harmonic representative, let's say $\omega+d\alpha$, for some $1$-form $\alpha$. Since $\omega$ is not exact (i.e. it is one of the generators of the second cohomology group of $\C P^2 \# \overline{\C P^2}$), there is a point $q$ where $\omega+d\alpha$ is non-zero. Then $\omega+d\alpha$ is non-trivial and in particular its rank $s$ at $q$ is $r_q\geq 2$. On the other hand, because $\omega^2=0$, we have $(\omega+d\alpha)^2=2 \omega \wedge d\alpha+d\alpha^2$ which is exact and thus it can not be a volume form by Stokes' Theorem. This means that there is a point $p$ where $(\omega+d\alpha)^2=0$ and then the rank of $(\omega+d\alpha)_p$ is $r_p\leq 2$. Therefore $\omega+d\alpha$ is non-trivial, $g$-harmonic, can not have constant rank $2$ and can not be symplectic. 
%\end{Example}

One can try to seek for more restricted ranks and ask the following:
\begin{Question}
Does $\C P^2 \# \overline{\C P^2}$ admit an intrinsically harmonic $2$-form of non-constant rank $r\leq 4$ which has rank $4$ at least in a point? 
\end{Question}

\begin{Question}
Does $\C P^2 \# \overline{\C P^2}$ admit an intrinsically harmonic $2$-form of non-constant rank $r$ such that $2\leq r\leq 4$ and there are points where $r=2$ and $r=4$? 
\end{Question}

Here $\C P^2 \# \overline{\C P^2}$ can be replaced by the blow-up of $\C P^2$ in $k$ point or, more generally, by a closed orientable symplectic $4$-manifolds which is non-minimal. 

Of course the more important is the following question, which should be compared with \cite[Conjecture $4.5$]{KA}: 
\begin{Question}
Characterize intrinsically generic (non-vanishing) harmonic $2$-forms on $4$-manifolds. The next step should be to drop the genericity condition. More generally what can one say in dimension greater than $4$?
\end{Question}

Coming back to $(n-1)$-forms, Remark \ref{rem: n-1 forms} suggests the following:
\begin{Question}
Is there an example of a nowhere-vanishing closed $(n-1)$-form, non-exact and non intrinsically harmonic?.
\end{Question}

Related to Corollary \ref{cor-simic} is the following question (see the final Question in \cite{SI}):
\begin{Question}
Is there some relation between two Riemannian metrics associated to two non-homologous cross-sections?
\end{Question}

Since symplectic forms are intrinsically harmonic with respect to a compatible metric, the following question is quite natural:
\begin{Question}
A symplectic form $\omega$ (in particular a K\"ahler form) is harmonic with respect to a compatible metric. Are there other interesting metrics making $\omega$ harmonic? Are two such metrics related in some sense?
\end{Question}

One could also check for maps {\it preserving} intrinsic harmonicity and therefore we have the:
\begin{Question}
Is it possible to give examples of a map between two manifolds such that pull back of intrinsically harmonic forms are intrinsically harmonic. Of course a diffeomorphism works just by pulling back the metric making a closed form harmonic.
\end{Question}

By \cite{K}, on a Riemannian manifold the fact that the wedge product of harmonic forms is harmonic, implies strong topological restrictions. Since intrinsic harmonicity does not fix the metric, one can ask what follows:
\begin{Question}
What one can say about manifolds for which the wedge product of intrinsically harmonic forms is still intrinsically harmonic?
\end{Question}

{\bf Acknowledgments}. The author wishes to express sincere gratitude to the organizers of the Second Conference on Differential Geometry, held in Fez in October 2024, for their warm hospitality and meticulous organization of the event. The author also thanks Giovanni Placini for his valuable feedback on an initial draft of this paper and the anonymous referee for his helpful comments and suggestions.


\begin{thebibliography}{10}

\bibitem{B}
G. Bande, {\sl Symplectic pairs and intrinsically harmonic forms}, Mathematics (2023), Volume 11, Issue 13, 2993.


\bibitem{BG}
G. Bande, P. Ghiggini, {\sl Holomorphic spheres and four-dimensional symplectic pairs},  J. Geom. Anal. {\bf 30} (2020), no. 1, 861--873.

\bibitem{BGK}
G. Bande, P. Ghiggini, D. Kotschick, {\sl A stability theorem
for contact and symplectic pairs}, Int. Math. Res. Not. (2004), no.
68, 3673--3688.

\bibitem{BK}
G. Bande, D. Kotschick, {\sl The Geometry of Symplectic pairs},
Trans.~Amer.~Math.~Soc.~\textbf{358} (2006), no. 4, 1643--1655.

%\bibitem{BK2}
%G. Bande, D. Kotschick, {\it The geometry of recursion operators}, Comm. Math. Phys. {\bf 280} (2008), no. 3, 737--749.
%
%\bibitem{BK3}
%G. Bande, D. Kotschick, {\it Contact pairs and locally conformally symplectic structures}, Harmonic maps and differential geometry, 85--98, Contemp. Math., 542, Amer. Math. Soc., Providence, RI, 2011. 

%\bibitem{Bru}
%M. Brunella, {\it Feuilletage holomorphes sur les surfaces complexes compactes}, Ann. Sci. \'Ecole Norm. Sup. (4) {\bf 30} (1997), no. 5, 569--594. 
%
%\bibitem{CC}
%A. Candel, L. Conlon, {\it Foliations I} Graduate Studies in Mathematics, {\bf 23}. American Mathematical Society, Providence, RI, 2000.


%\bibitem{DC}
%M. P. do Carmo, {\sl Riemannian geometry}, Mathematics: Theory \& Applications. Birkh\"auser Boston, Inc., Boston, MA, 1992.
%


%\bibitem{HLS}
%H.~Hofer, V.~Lizan and J.-C.~Sikorav, \textit{On genericity for holomorphic curves in four-dimensional almost-complex manifolds},  J. Geom. Anal. {\bf 7} (1997), no. 1, 149--159. 
%

\bibitem{C}
E. Calabi, {\sl An intrinsic characterization of harmonic one-forms}, Global Analysis (Papers in Honor of K. Kodaira) pp. 101--117, Tokyo, 1969 

\bibitem{FKL} M. Farber, G. Katz, J. Levine, Jerome, {\sl Morse theory of harmonic forms}, Topology {\bf 37} (1998), 469--483.


\bibitem{F}
E. Fran\c ca, {\sl Intrinsically harmonic forms and characterization of flat
circle bundles}, PhD Thesis, 2020.

\bibitem{EF}
E. Fran\c ca, D. Finamore, {\sl A characterization of the n-dimensional torus via intrinsically harmonic forms}, arXiv:2205.15807, 2023.

\bibitem{EM}
E. Fran\c ca, F. Mercuri, {\sl Intrinsically harmonic forms and flat
circle bundles}, arXiv:2205.14190, 2022.


\bibitem{G}
C.~Godbillon, {\sl Feuilletages}, Birkh\"auser Verlag 1991.

\bibitem{H}
K. Honda,  {\sl On harmonic forms for generic metrics}, PhD Thesis, Princeton University, Princeton, 1997

\bibitem{KA} G. Katz, {\sl Harmonic forms and near-minimal singular foliations}, Comment. Math. Helv. {\bf 77} (2002), 39--77.

\bibitem{K}
D.~Kotschick, {\sl On products of harmonic forms}, Duke Math. J. {\bf 107} (2001), 521--531.

%\bibitem{KM}
%D.~Kotschick, S.~Morita, {\it Signatures of foliated surface
%bundles and the symplectomorphism groups of surfaces}, Topology
%{\bf 44} (2005), no. 1, 131--149.

%\bibitem{MC}
%D. McDuff, {\it The structure of Rational and Ruled Symplectic $4$-manifolds},  J. Amer. Math. Soc. {\bf 3} (1990), no. 3, 679--712.
%
%\bibitem{MC2}
%D. McDuff, {\it The local behaviour of holomorphic curves in almost- complex $4$-manifolds}, J. Diff. Geom. {\bf 34} (1991) 211--358.

\bibitem{MS}
D. McDuff, D. Salamon, {\sl Introduction to Symplectic Topology. Second edition.} Oxford Mathematical Monographs. The Clarendon Press, Oxford University Press, New York, 1998.

%\bibitem{MCS2}
%D. McDuff, D. Salamon, {\sl J-holomorphic curves and symplectic topology. Second edition.} American Mathematical Society Colloquium Publications, 52. American Mathematical Society, Providence, RI, 2012. 
%
%\bibitem{MP}
%V. Mu\~noz, F. Presas, {\it Symplectic Foliations in Four-Manifolds}, Proceedings of the X Fall Workshop on Geometry and Physics.

%\bibitem{R}
%H. Rummler, {\it Quelques notions simples en g\'eom\'etrie riemannienne et leurs applications aux feuilletages compacts}, Comment. Math. Helv. {\bf 54} (1979) 224--239.
%
%\bibitem{S}
%A. Scorpan, {\it Existence of foliations on 4-manifold}, Algebr. Geom. Topol. {\bf 3} (2003), 1225--1256.
%
%\bibitem{W2}
%C.~Wendl, {\sl Lectures on Holomorphic Curves in Symplectic and Contact Geometry}, https://arxiv.org/abs/1011.1690, 2014.
%
%\bibitem{W}
%C.~Wendl, {\sl Holomorphic Curves in Low Dimensions}, https://www.mathematik.hu-berlin.de/~wendl/pub/rationalRuled.pdf, 2017.
\bibitem{SI} S. Simi\'c, {\sl Cross-sections to flows and intrinsically harmonic forms}, Dyn. Syst. {\bf 38} (2023), 314--319.

\bibitem{V} E. Volkov, {\sl Characterization of intrinsically harmonic forms}, J. Topol. {\bf 1} (2008), no. 3, 643--650. 

\end{thebibliography}
\end{document}